\documentclass{article}

\usepackage{amssymb,amsmath,eucal}
\usepackage[dvips]{graphicx}

\begin{document}

\def\Sp{\mathrm {Sp}}
\def\U{\mathrm U}
\def\SOS{\mathrm {SO}^*}

\def\R{{\mathbb R}}
\def\C{{\mathbb C}}
\def\T{{\mathbb T}}

\def\M{\mathcal M}

\def\B{{\rm B}}
\def\W{{\rm W}}
\def\ov{\overline}
\def\wt{\widetilde}

\def\Symp{\mathrm {Symp}}
\def\SSymp{\mathrm {SSymp}}
\def\Map{\mathrm {Map}}
\def\Ams{\mathrm {Ams}}

\def\fB{\mathfrak B}

\def\cF{\mathcal F}

\def\fT{\mathfrak T}
\def\fU{\mathfrak U}

\def\phi{\varphi}
\def\kappa{\varkappa}

\def\cL{\mathcal L}

\def\epsilon{\varepsilon}

\def\konets{{\hfill$\square$}}

\def\Ang{\mathrm{Ang}}
\def\ang{\mathrm{ang}}
\def\Aut{\mathrm{Aut}}

\def\tr{{\rm tr\,\,}}
\def\Re{{\rm Re\,\,}}
\def\Im{{\rm Im\,\,}}

\def\Z{{\mathbb Z}}

\def\BB{{\cal B}}

\def\SL{{\rm SL}}

\def\Move{\mathrm{Move}}


\renewcommand{\theequation}{\arabic{section}.\arabic{equation}}
\newcounter{punct}[section]


\renewcommand{\thepunct}{\thesection.\arabic{punct}}
\def\punct{\refstepcounter{punct}{\arabic{section}.\arabic{punct}.  }}


\def\SS{\smallskip}

\newtheorem{theorem}{Theorem}[section]
\newtheorem{proposition}[theorem]{Proposition}
\newtheorem{prop}[theorem]{Proposition}
\newtheorem{lemma}[theorem]{Lemma}
\newtheorem{cor}[theorem]{Corollary}
\newtheorem{corollary}[theorem]{Corollary}
\newtheorem{observation}[theorem]{Observation}

\begin{center}

{\Large \bf Central extensions
of  groups of symplectomorphisms}

\medskip

\large\sc Yuri A. Neretin%
\footnote{Partially supported by the grant NWO.047.017.015  and
the grant FWF, project P19064}

\end{center}

\bigskip

{\small We construct  canonically
defined central extensions of
 groups of symplectomorphisms.
We show that this central extension is nontrivial in the case
of a torus of dimension $\ge 6$ and in the case of
a two-dimensional surface of genus $\ge 3$.}

\setcounter{equation}{0}

\bigskip

\section{Formulation of results}

\medskip

Central extensions of the groups of symplectomorphisms
discussed in this paper appeared as a byproduct
in \cite{Ner2}.  Here we prove several
nontrivality and triviality theorems concerning
this cocycle.

\smallskip


{\bf\punct\label{ss1.1} Preliminaries.
 Cocycle on the symplectic group $\Sp(2n,\R)$.}
We define the real symplectic
group $\Sp(2n,\R)$ as the group of real matrices
\begin{equation}
g=\begin{pmatrix} A&B\\C&D\end{pmatrix}
\label{eq1.1}
\end{equation}
preserving the standard skew-symmetric bilinear form
$K:=\begin{pmatrix} 0 &1 \\ -1&0\end{pmatrix}$,
i.e.
\begin{equation}
g^t K g=K
\label{eq1.2}
\end{equation}
The complex symplectic group
$\Sp(2n,\C)$
is the group
of complex matrices satisfying  the same condition
 (\ref{eq1.2}).

Consider the block $(n+n)\times (n+n)$ matrix $J\in\Sp(2n,\C)$
given by
\begin{equation}
J=\frac1{\sqrt 2}
\begin{pmatrix} 1&i\\i&1
\end{pmatrix}
\label{eq1.3}
\end{equation}

For $g\in \Sp(2n,\R)$, we consider
the matrix $J^{-1}gJ\in\Sp(2n,\C)$,
this matrix has the structure
$$J^{-1}gJ=
\begin{pmatrix}
\Phi&\Psi\\ \ov\Psi &\ov\Phi
\end{pmatrix}
$$
where bar means the element-wise complex conjugation.
We denote
$$\Phi=\Phi(g);\qquad \Psi=\Psi(g)$$
%

We define the {\it Berezin cocycle}
$$
c:\Sp(2n,\R)\times\Sp(2n,\R)\to\R
$$
by%
\footnote{In conciderations of Subsection \ref{ss4.2},
 this formula
arises in a natural way} 
\begin{equation}
c(g_1,g_2)=\Im
\tr\ln\Bigl[
\Phi(g_1)^{-1}\Phi(g_1g_2)(\Phi(g_2)^{-1}\Bigr]
\label{eq1.4}
\end{equation}

Below (Theorem \ref{th2.1}) we show that the matrix  in the brackets
has the form $1+Z$, where $\|Z\|<1$. Then
$\ln(1+Z):=Z-Z^2/2+Z^3/3-\dots$ and hence our expression
is well defined.

\smallskip

The cocycle $c$ defines a central extension of $\Sp(2n,\R)$.
In other words, the set $\Sp(2n,\R)\times \R$ with the multiplication
$$(g,x)\cdot (h,y)=(gh,x+y+c(g,h))$$
is a group.

\smallskip


{\bf\punct\label{ss1.2} Groups of symplectomorphisms. Notation.}
 Consider  a $2n$-dimensional
symplectic manifold $M$.
Define the following groups

--- $\Symp(M)$ is  the group of
all $C^\infty$-smooth {\it compactly supported}
 symplectomorphisms
of $M$. If the manifold $M$ itself is compact,
then $\Symp(M)$ is the group of all symplectomorphisms
of $M$.

--- $\SSymp(M)$ is the connected component of
$\Symp(M)$ containing unit $e$.

--- $\SSymp^\sim(M)$ is the {\it universal covering group}
of $\SSymp(M)$.

--- More generally, we denote
the universal covering of any connceted group $G$ by $G^\sim$

--- $\Map(M)$ is  the {\it mapping class
group} $\Symp(M)/\SSymp(M)$.

\smallskip

We have a natural topology on $\Symp(M)$ and hence we have a natural
Borel structure on $\Symp(M)$. In particular, we have a notion
of a {\it measurable function} on $\Symp(M)$. Let $F$ 
be a measurable function on $\Symp(M)$, let $U\subset \R^N$ be an open
domain and $\psi:U\mapsto \Symp(M)$ be a smooth map.
Then $F\circ\psi$ is a measurable function on $U$ in the usual
sense.   

\smallskip 


{\bf\punct\label{ss1.3} Central extension
of the group of symplectomorphisms.}
Equip the space $\R^{2n}$  with
the standard symplectic structure.
Consider an open set $\Omega\subset\R^{2n}$
and a symplectic embedding $\iota:\Omega\to M$
such that the measure of $M\setminus \iota(\Omega)$
is zero.

\smallskip

{\sc Remarks.} a) We admit disconnected sets $\Omega$.

b) It is pleasant (but not necessary) to think that
$M\setminus \iota(\Omega)$ is a union
of submanifolds.

\smallskip

Any element
$g\in\Symp(M)$ induces a
transformation $\iota^{-1}g\iota$ of $\Omega$
defined almost everywhere.
Denote the group of all such transformations
by $\Symp(M,\Omega,\iota)$.
By the definition, $\Symp(M,\Omega,\iota)\simeq\Symp(M)$.
This group contains $\Symp(\Omega)$ as a proper
subgroup.

 For $q\in\Symp(M,\Omega,\iota)$
and $x\in\Omega$,
we denote by $q'(x)$ its Jacobi matrix
at the point $x$.
We define the 2-cocycle $C(q_1, q_2)$,
where
$q_1$, $q_2\in\Symp(M,\Omega,\iota)$, by the formula%
\footnote{This cocycle is induced from a cocycle on
the group of measurable currents, see below
 Subsection \ref{ss2.4}.
Also, below we propose a coordinate-less description
of the cocycle $C$, see Subsection \ref{ss2.6}.}
\begin{multline}
C(q_1, q_2)=\\=
\int_\Omega \Im\tr\ln\Bigl\{
\Phi^{-1}\bigl[q'_1(q_2(m)\bigr]\,\,
\Phi\bigl[q'_1(q_2(m))q'_2(m)\bigr]\,\,
\Phi^{-1}\bigl[q_2'(m)\bigr]\Bigr\}
dm=\\
\int_\Omega c\bigl(q_1'(q_2(m)), q_2'(m)\bigr)\,dm
\label{eq1.5}
\end{multline}

\begin{theorem}
\label{th1.1}
a) The expression $C(q_1, q_2)$ defines an element of
the second cohomology group $H^2(\Symp(M),\R)$.
In another words, the space $\Symp(M)\times \R$
with the product
\begin{equation}
(q_1,x_1)\times (q_2,x_2)=(q_1\circ q_2,x_1+x_2+C(q_1,q_2))
\label{eq1.6}
\end{equation}
is a group.

b)  This central extension of $\Symp(M)$ does not
depend on a choice of the domain $\Omega$
and the map $\iota$.
\end{theorem}


\smallskip

{\bf\punct\label{ss1.4} Triviality results for the cocycle $C$.}
 We have a map from the universal covering group
  $\SSymp^\sim(M)$ to $\SSymp(M)$
and hence we can consider the cocycle $C$ as a
cocycle on $\SSymp^\sim(M)$.

\smallskip

\begin{proposition}
\label{prop1.2}
  Let $\Xi\subset \R^{2n}$ be an open domain.
Then the central extension of $\SSymp^\sim(\Xi)$
defined by the cocycle $C$ is trivial.
\end{proposition}

\smallskip

Let $M$ be a symplectic manifold.
 Consider an almost complex structure
  on the tangent bundle of $M$
 compatible with the symplectic structure
 (in particular, we obtain $n$-dimensional
  complex vector bundle).
 Assume that the corresponding Hermitian metric
 is positive definite.

For the complex bundle obtained in this way,
consider its $n$-th exterior power $L$.

\smallskip

\begin{proposition}
\label{prop1.3} 
 If the
complex line bundle
$L$ on $M$ is trivial,
then the central extension of $\SSymp^\sim(M)$
defined by the cocycle $C$ is trivial.
\end{proposition}

\begin{corollary}
\label{cor:add1}
For
a noncompact 2-dimensional surfaces $\M$ 
of a finite genus,
 our cocycle is trivial  on
$\SSymp(\M)$.
\end{corollary}

Indeed, in this case, the group $\SSymp(\M)$ is contractible,
 see \cite{EE}, hence $\SSymp^\sim(\M)=\SSymp(\M)$.
 Also, in this case, 
the line bundle $L$ is trivial. Thus,
 by Proposition \ref{prop1.3}, our central extension
is trivial.

\smallskip


{\bf \punct\label{ss1.5} Nontriviality results.}

\nopagebreak

\smallskip

\begin{theorem}
\label{th1.4} 
 For a two-dimensional
 oriented (compact or noncompact) surface $\M_g$ 
  of genus $g\ge 3$,
 our central extension of $\Symp(\M_g)$
 is nontrivial in measurable cohomologies.%
\end{theorem}

 \smallskip

Further, consider the standard lattice $\Z^{2n}$ in the standard
symplectic space $\R^{2n}$. Consider the torus
$\T^{2n}:=\R^{2n}/\Z^{2n}$.
Denote by $\Sp(2n,\Z)$ the subgroup in
$\Sp(2n,\R)$ consisting
of all matrices (\ref{eq1.1}) with integer elements.

The action of $\Sp(2n,\Z)$ on $\R^{2n}$ induces
the symplectic action of $\Sp(2n,\Z)$ on $\T^{2n}$.

\smallskip

\begin{theorem} 
\label{th1.5}
 The central extension of $\Symp(\T^{2n})$
defined by the cocycle $C$ is nontrivial
for $n\ge 3$.
\end{theorem}

\smallskip

In fact, we prove that this extension is nontrivial
on the countable subgroup $\Sp(2n,\Z)\subset\Symp(\T^{2n})$.
The latter statement is a kind of a rigidity
theorem for  lattices.

 \smallskip


{\bf \punct Some discussion.
\label{ss1.6}}

a) Central extensions of  of $\SSymp(M)$
were discussed in several works, see
Kostant \cite{Kos}, Brylinski \cite{Bry},
Ismagilov \cite{Ism0}, \cite{Ism}, Haller, Vizman \cite{HV}.
 Our construction differs from
these constructions.

b) Consider a surface $\M_g$ of genus $\ge 3$.
 It seems that our central extension
in this case must be related to
the Harer central extension \cite{Har}
of the mapping class group $\Map(\M_g)$,
see also \cite{Ger}.

But this relation now is not clear.
Indeed, $\Symp(\M_g)$ is not a semidirect product
$\Map(\M_g)\ltimes \SSymp(\M_g)$
and hence our construction does not induce
automatically an extension of $\Map(\M_g)$.

Second, the central extension of $\Map(\M_g)$
induces a central extension of $\Symp(\M_g)$.
Unfortunately, I do not know,
gives our construction the same result or not.


c) Symplectic mapping class groups
in higher dimensions were discussed by Seidel, see \cite{Sei};
for results and references on topology
of groups of symplectomorphisms, see
 surveys \cite{McD}, \cite{McD2}.

d) For a two-dimensional surface $\M_g$,
our central extension
can be realized in a unitary representation of $\Symp(\M_g)$,
this construction was obtained in \cite{Ner1}.
I do not believe that a realization in a unutary representation
  is possible
for dimensions $\ge 4$ (realizations in  nonunitary
 representations exist).

\smallskip


{\bf\punct Structure of the paper.
\label{ss1.7}}
Details of construction of the cocycle $C$
and proofs of Theorem \ref{th1.1} and 
Propositions \ref{prop1.2}--\ref{prop1.3}
are contained in Section 2.

Theorems \ref{th1.4} and \ref{th1.5}
 are proved in Sections 3 and 4
respectively.

\smallskip

{\bf Acknowledgments.}
I am very grateful to V.Fock who explained me how to
translate the measure-theoretic construction of \cite{Ner2}
to the language of vector bundles (see Subsection \ref{ss2.6}).
I thank R.S.Ismagilov who  simplified
the proof of Theorem 1.4 respectively 
the preprint variant of the work.

The main part of this work was  done during my visits
to Lyon and Grenoble in December 1999 and
to Vienna in December 2003. I thank C.Roger,
V.Sergiescu, and P.Michor
 for discussions and hospitality.
I also thank N.V.Ivanov and S.Haller for discussions
and references.

\medskip


\section{Constructions of cocycles}

\medskip

\nopagebreak

\setcounter{equation}{0}


{\bf\punct Central extensions. Preliminaries.\label{ss2.1}}
{\bf A.}
 Let $G$ be a group,
let $A$ be an Abelian group (in our work,
$A$ is the additive group of $\R$, $\Z$, or
the  circle $\T:=\R/\Z$, we write the operation
in $A$ in the additive form).
A {\it central extension}
of $G$ by the group $A$
(or $A$-extension of $G$) is a group $\wt G$
such that $A$ is a central subgroup in $\wt G$
and $\wt G/A\simeq G$.

The set $\wt G$ can be identified noncanonically with
the product $G\times A$ and the homomorphism
 $\wt G\to G$ can be identified with the projection
 $G\times A\to G$. Then the multiplication in $G\times A$
 must have the form
\begin{equation}
 (g_1, a_1)\cdot(g_2,a_2)=
 \bigl(g_1g_2, a_1+a_2+c(g_1,g_2)\bigr)
\label{eq2.1} 
\end{equation}
 where the function $c:G\times G\to A$ (a {\it 2-cocycle})
 satisfies the identities
 \begin{align}
 c(e,g)=c(g,e)=0\qquad
    \text{for all $g$ }
\label{eq2.2}
;\\
c(g_1,g_2)+c(g_1g_2,g_3)=
c(g_1,g_2g_3)+c(g_2,g_3)
\label{eq2.3}
\end{align}
 where
  $e$ is the unit of $G$.
  The first condition means that $(e,0)$ is the unit
  of $\wt G$.
The second condition is equivalent
to the associativity of the multiplication
(\ref{eq2.1}).

{\bf B.}
If we change the identification $G\times A$ and $\wt G$,
then $c(g_1,g_2)$ is changed according  the rule
\begin{equation}
c(g_1,g_2)\mapsto c(g_1,g_2)-
\gamma(g_1)-\gamma(g_2)+\gamma(g_1g_2)
\label{eq2.4}
\end{equation}
where $\gamma:G\to A$ is some function
such that $\gamma(e)=0$.

The central extension is { \it trivial}
if $c(g_1,g_2)$ can be transformed
to 0 by the operation (\ref{eq2.4}), i.e.,
\begin{equation}
c(g_1,g_2)=\gamma(g_1)+\gamma(g_2)-\gamma(g_1 g_2)
\label{eq2.5}
\end{equation}
In this case $\wt G=G\times A$ and we say that
$\gamma(g)$ is a {\it trivializer} of $c$.

{\bf C.}
If $\gamma$, $\gamma'$ are two trivializers
of the same cocycle $c$, then
$\gamma-\gamma'$ is a homomorphism $G\to A$.

{\bf D.}
The additive group of functions $c(g_1,g_2)$
satisfying (\ref{eq2.2})--(\ref{eq2.3}) is denoted by $C^2(G,A)$
({\it group of cocycles}); the group of functions
having the form (\ref{eq2.5}) is denoted by $B^2(G,A)$
({\it the group of coboundaries}).
The {\it second cohomology group} is the factor-group
$$
H^2(G,A):=C^2(G,A)/B^2(G,A)
$$

{\bf E.}
Let $B$ be an Abelian group, $A$ be its subgroup,
$C=B/A$ be the factor-group. Then we have the obvious maps
\begin{equation}
C^2(G,A)\to C^2(G,B)\to C^2(G,C);
\label{eq2.6}
\end{equation}
The first map
means that an $A$-valued function $c$ is also a $B$-valued function;
considering a composition of a function
$G\times G\to B$ and the homomorphism $B\to C$,
we obtain the second map.
We also have the corresponding map in
cohomologies
\begin{equation}
H^2(G,A)\to H^2(G,B)\to H^2(G,C)
\label{eq2.7}
\end{equation}

{\bf F.}
Let $G$, $G'$ be groups,
let $\theta:G'\to G$ be a homomorphism.
Then $\theta$ induces a natural map
$C^2(G,A)\to C^2(G',A)$, i.e., for a
cocycle $c\in C^2(G,A)$ we consider
the cocycle $c(\theta(g_1'),\theta(g_2'))\in C^2(G',A)$.
Hence we also have a map of cohomologies
$$H^2(G,A)\to H^2(G',A)$$

{\bf G.} Let $G$ be a group, fix $h\in G$.

The cocycle $c(h^{-1}g_1 h,h^{-1}g_2 h)$
is equivalent to $c(g_1,g_2)$, see \cite{Bro}, III.8.1.

\smallskip

{\bf H.} Now let $G$, $A$ be  topological groups. We say that a cocycle
$c\in H^2(G,A)$ is nontrivial in measurable cohomologies
if it can not be trivialized 
(see (\ref{eq2.6})) by a measurable trivializer $\gamma$.

\smallskip


{\bf\punct A model of $\Sp(2n,\R)$.\label{ss2.2}}
In \ref{ss1.1}, we  realized the group
$\Sp(2n,\R)$
as the group
 of complex  matrices
having the block structure
\begin{equation}
g=\begin{pmatrix}
\Phi&\Psi\\ \ov\Psi &\ov\Phi
\end{pmatrix}
\label{eq2.8}
\end{equation}
preserving the
skew-symmetric bilinear form  $K$
with the matrix
$\begin{pmatrix} 0 &1 \\ -1&0\end{pmatrix}$.

 Then the matrix (\ref{eq2.8}) also

--- preserves the indefinite Hermitian form  $M$
with the matrix
$\begin{pmatrix}1&0\\0&-1\end{pmatrix}$,
i.e.,
\begin{equation}
g\begin{pmatrix} 0 &1 \\ -1&0\end{pmatrix}g^*=
g^*\begin{pmatrix} 0 &1 \\ -1&0\end{pmatrix}g=
\begin{pmatrix} 0 &1 \\ -1&0\end{pmatrix}
\label{eq2.9}
\end{equation}

-- preserves the real subspace $V\subset \C^n$,
 consisting of the vectors $(h,\ov h)$,
 moreover,  matrices (\ref{eq2.8}) commute with the
antilinear operator
 $(p,q) \mapsto (\ov q,\ov p)$.

 \medskip

Now let us give a
 coordinateless description
of this realization of $\Sp(2n,\R)$.
Consider an $n$-dimensional complex space $V$
equipped with a positive definite Hermitian form
$H(\cdot,\cdot)$. Denote by $V_\R$ the same space
considered as $2n$-dimensional
 linear space over $\R$.
The operator $v\mapsto iv$ in $V$ is
also a linear operator in $V_\R$, we denote it
by $I$.

Denote by $\{\cdot,\cdot\}$
the imaginary part
of the Hermitian form  $H$,
 it is a skew-symmetric
bilinear form on $V_\R$.
The group preserving this form is $\Sp(2n,\R)$.

Consider the complexification
$(V_\R)_\C$ of the space $V_\R$. It is a
$2n$-dimensional
complex linear space equipped with
several additional structures

1) We have an operator
$I$ such that $I^2=-1$.

2) Since $I^2=-1$, the eigenvalues of $I$ are $\pm i$.
Denote by $V_\pm$ the corresponding eigenspaces.
Then $V=V_+\oplus V_-$.

3) We have
the operation $Q$ of complex conjugation
$Q:v+iw\mapsto v-iw$, where $v,w\in V_\R$.
It satisfies $QV_\pm=V_\mp$

Now we consider the action of $\Sp(2n,\R)$
in $(V_\R)_\C$.
 An operator $g\in\Sp(2n,\R)$
 preserves
the bilinear form in $(V_\R)_\C$ and commutes
with the complex conjugation.

Representing it as a block operator
$V_+\oplus V_-\to V_+\oplus V_-$,
we obtain
the block  matrix  representation (\ref{eq2.8})

The Hermitian  form $M$ on $(V_\R)_\C$ is
$$M(v+iw,v'+iw')=\{v,w'\}- \{w,v'\}+ i\{v,v'\}+i\{w,w'\}$$
it is more natural to say that we extend $i\{\cdot,\cdot\}$
as an Hermitian form from $V_\R$ to $(V_\R)_\C$.



\SS

{\bf \punct The Berezin cocycle on $\Sp(2n,\R)$.\label{ss2.3}}

\smallskip

\begin{theorem}
\label{th2.1}
 a)  The function
$\Sp(2n,\R)\times\Sp(2n,\R)\to \R$ given by
\begin{equation}
c(g_1,g_2)=\Im
\tr\ln\Bigl[
\Phi(g_1)^{-1}\Phi(g_1g_2)(\Phi(g_2)^{-1}\Bigr]
\label{eq2.10}
\end{equation}
 is well-defined.

 b)  The function $c(g_1,g_2)$ is a 2-cocycle.

 c)  The $\R$-valued cocycle $\frac 1{2\pi}c(g_1,g_2)$
 can be reduced to
 a $\Z$-valued cocycle.%
 \footnote{I.e., $c$ is contained in the image of the map
 $H^2(\Sp(2n,\R),\Z)\to H^2(\Sp(2n,\R),\R)$.}
 The corresponding $\Z$-extension of $\Sp(2n,\R)$
  coincides with
 the universal covering group of $\Sp(2n,\R)$.
\end{theorem}

  d) {\it The cocycle $c$ is uniformly bounded
on $\Sp(2n,\R)\times\Sp(2n,\R)$ }
 \begin{equation}
 |c(g_1,g_2)|<n\pi/2
\label{eq2.11}
\end{equation}

 \smallskip

 {\sc Remark.} In \cite{Ber1}, Berezin wrote the
 following $\T$-cocycle on $\Sp(2n,\R)$ and on its
 infinite-dimensional analogue
 $$\sigma(g_1,g_2)=\det (1+\Phi(g_1)^{-1}\Psi(g_1)\cdot\ov\Psi(g_2)\Phi(g_2)^{-
1})^{-1/2}
$$
Trivial  calculation shows that
$$\sigma(g_1,g_2)=\exp\{-c(g_1,g_2)/2\}$$
This cocycle can be trivialized on the two-sheet covering
of $\Sp(2n,\R)$.
  For $n=1$, an explicit formula for
 $\R$-valued cocycle $c$
 was written by Guichardet \cite{Gui}.

 \smallskip

 {\sc Proof.}
 a) We have
 $$\Phi(g_1 g_2)=\Phi(g_1)\Phi(g_2)+\Psi(g_1)\ov\Psi(g_2)$$
 Hence
 \begin{equation}
 \Phi(g_1)^{-1} \Phi(g_1g_2)\Phi(g_2)^{-1}=
 1+\Phi(g_1)^{-1}\Psi(g_1)\cdot\ov\Psi(g_2)\Phi(g_2)^{-1}
\label{eq2.12}
 \end{equation}
 Equations (\ref{eq2.9}) imply
 \begin{equation}
 \Phi\Phi^*-\Psi\Psi^*=1;\qquad
 \Phi^*\Phi-\Psi^t\ov\Psi=1
\label{eq2.13}
 \end{equation}
 Hence $\Phi$ is invertible and
 we have the following inequalities for norms
\begin{equation}
\|\Phi^{-1}\Psi\|<1,\qquad
\|\ov \Psi\Phi^{-1}\|<1
\label{eq2.14}
\end{equation}
Thus,
$$
\|\Phi(g_1)^{-1}\Psi(g_1)\cdot\ov\Psi(g_2)\Phi(g_2)^{-1}\|<1
$$
We define the logarithm of (\ref{eq2.12}) by
$$
\ln(1+Z):=\sum_{n=1}^\infty (-1)^{n+1} Z^n/n
$$
Since the norm of our $Z$ is $< 1$, our series converges.

b) Let $A$, $B$ be $n\times n$ matrices and
$\|A-1\|$, $\|B-1\|$ be sufficiently small.
Then
$$
\tr\ln(AB)=\tr\ln A+\tr\ln B
$$
Hence, for $g_1$, $g_2$ lying  in a small neighborhood
of the unit we can write
$$
\tr\ln\Bigl[
\Phi(g_1)^{-1}\Phi(g_1g_2)\Phi(g_2)^{-1}\Bigr]
=
-\tr\ln\Phi(g_1)+\tr\ln\Phi(g_1g_2)-\tr\ln\Phi(g_2)
$$
Now the cocycle identity (\ref{eq2.3}) became
trivial for $g_1$, $g_2$, $g_3$
lying in a sufficiently small
neighborhood of unit.

But all our expressions
are real analytic and the group $\Sp(2n,\R)$
is connected. Hence the cocycle  identity (\ref{eq2.3})
is valid for
all $g_1$, $g_2$, $g_3\in \Sp(2n,\R)$.

A verifying of (\ref{eq2.2}) is trivial.

\smallskip

c) As we have seen $\det\Phi\ne 0$. Assume
\begin{equation}
\gamma(g):=\Im\tr\ln \Phi(g)=\arg\ln \det\Phi(g)
\label{eq2.15}
\end{equation}
where  $0\le \arg z< 2\pi$.
Obviously,
$$\alpha:=\frac 1{2\pi}\Bigl(
c(g_1,g_2)+\gamma(g_1)+\gamma(g_2)-\gamma(g_1g_2)
\Bigr)\in \Z$$
since $\exp(2\pi \alpha)=1$,
and the first statement is proved.

The $\Z$-central extension of
 $\Sp(2n,\R)$ obtained in this way
can be considered as the subset $\Sp(2n,\R)^\sim
\subset
\Sp(2n,\R)\times \R$ consisting of pairs
$$(g,x), \text{\qquad where
$g\in\Sp(2n,\R)$, $x\in\R$ and $\det \Phi(g)/|\det\Phi(g)|=e^{ix}$}
,$$
 the multiplication is given by (\ref{eq2.1}).
Obviously, the projection
$\Sp(2n,\R)^\sim\to\Sp(2n,\R)$ is a covering map.
For $y\in \R$ consider the matrix $g(y)$ with
$\Psi(y)=0$ and $\Phi$ being the diagonal matrix
with entries
$e^{iy}$, $1$,\dots, $1$. The map
$y\mapsto (g(y),y)$ is a continuous map
$\R\to\Sp(2n,\R)^\sim$ and hence $\Sp(2n,\R)^\sim$
is connected.
Thus $\Sp(2n,\R)^\sim$ is a covering group for $\Sp(2n,\R)$.

It remains to notice that
the unitary group $\U(n)$ is a deformation retract
of $\Sp(2n,\R)$ and the loop $g(y)$, $y\in[0,2\pi]$,
is a generator of the fundamental group
$\pi_1(\U(n))$. Hence $\Sp(2n,\R)^\sim$ is a universal covering
of $\Sp(2n,\R)$.

\smallskip

d) Let $Z$ be a matrix satisfying
 $\|Z\|<1$. Let $\lambda_j(Z)$ be its eigenvalues.
Obviously, $|\lambda_j|<1$.
Let us prove that
\begin{equation}
\tr\ln (1+Z)=\sum \ln(1+\lambda_j(Z))
\label{eq2.16}
\end{equation}
For a self-adjoint matrix $Z$ this identity is obvious.
The both sides of (\ref{eq2.16}) 
are complex analytic in $Z$
and hence this identity is valid if $\|Z\|<1$.

But $|\lambda_j(Z)|<1$ implies $|\Im(1+\lambda_j)|<\pi/2$
and we obtain (\ref{eq2.11}).
\hfill $\square$

\smallskip

{\sc Remark.} Obviously,
the expression
$$
\Re\tr\ln\Bigl[
\Phi(g_1)^{-1}\Phi(g_1g_2)(\Phi(g_2)^{-1}\Bigr]
$$
also satisfies the cocycle equation (\ref{eq2.3}),
 but this cocycle is trivial,
since its trivializer
$$
\gamma(g)=\Re\tr\ln\Phi(g)=
\ln|\det(\Phi(g))|
$$
is well defined (but similar expression (\ref{eq2.15}) 
is multivalued).

\smallskip


{\bf\punct Groups $\fB(X,G)$.\label{ss2.4}} Denote by $X$
some Lebesgue space with a continuous measure%
\footnote{In fact, any such space is equivalent to
a segment $[a,b]\subset \R$ equipped with
the Lebesgue measure or to the whole
line $\R$.} 
 $\mu$.
For example, we can consider $X$ being an arbitrary
symplectic manifold.
Denote by $\Ams(X)$ the group of all measure-preserving
maps from the space $X$ to itself
(by definition, these maps are defined almost everywhere).

Let $G$ be an arbitrary group (below $G=\Sp(2n,\R)$).
Denote by $\cF(X,G)$
the group of all measurable functions
$f:X\to G$ such that
$$
f(x)=1\qquad \text{outside a set of finite measure}
$$
If $X$ itself has a finite measure, then we
can forget the last condition.

Consider the semidirect product
$$
\fB(X,G)=\Ams(X)\ltimes\cF(X,G)
$$
Elements of the group $\fB(X,G)$ are pairs
$\{p,h\}\in\Ams\times \cF(X,G)$
and the product is given by
$$
\{p_1(x),h_1(x)\} * \{p_2(x),h_2(x)\}:=
\{p_1(p_2(x)),h_1(p_2(x))h_2(x)\}
$$

{\sc Remark.} Let the group $G$ acts on a manifold $Y$
by transformations $y\mapsto yg$. 
Consider the space $\cL(X,Y)$ of all $Y$-valued
measurable functions on $X$. The group $\fB(X,G)$
acts in this space by the transformations
\begin{equation}
T(p,h)f(x)=f(p(x))h(x)
\label{eq2.17}
\end{equation}

\smallskip


{\bf\punct \label{ss:ext-fB} 
A central extension of $\fB(X,G)$.}
Let $c\in H^2(G,\R)$ be a bounded cocycle.
We define the function $C:\fB(X,G)\times\fB(X,G)\to\R$
by
\begin{equation}
C(\{p_1,h_1\},\{p_2,h_2\}=
\int_{X}
c\bigl(h_1\circ p_2(x), h_2(x)\bigr)\,d\mu(x)
\label{eq2.18}
\end{equation}

\smallskip

\begin{theorem}
\label{th2.2}
 The function $C$ is a 2-cocycle
on $\fB(X,G)$.
\end{theorem}

We denote the corresponding central extension
of $\fB(X,G)$ by $\wt\fB(X,G)$.


\smallskip

{\sc Proof.} We directly verify the 
cocycle identity (\ref{eq2.3}).
First,
\begin{gather}
C\Bigl(\{p_1,h_1\}\,,\,\{p_2,h_2\}\Bigr)+
C\Bigl(\{p_1,h_1\}*\{p_2,h_2\}\,,\,\{p_3,h_3\}\Bigr)=
\label{eq2.19}\\
=\int_X
c\Bigl([h_1\circ p_2](x)\,\,  ,   \,\, h_2(x)\Bigr)\,d\mu(x)
+
\label{eq2.20}
\\+
\int_X
c\Bigl([h_1\circ p_2\circ p_3](x)\cdot [h_2\circ p_3](x)
\,\, ,
\,\,
     h_3(x)\Bigr)\,d\mu(x)\nonumber
\end{gather}

Secondly,
\begin{gather}
C\Bigl(\{p_2\,,\,h_2\}\,,\,\{p_3\,,\,h_3\}\Bigr)
+
C\Bigl(\{p_1,h_1\}\,,\,\{p_2,h_2\}*\{p_3,h_3\}\Bigr)
=
\label{eq2.21}
\\=
\int_X
c\Bigl([h_2\circ p_3](x)\,\,,\,\, h_3(x)\Bigr)\,d\mu(x)
+
\label{eq2.22}
\\+
\int_X
c\Bigl([h_1\circ p_2\circ p_3](x)\,\,,\,\,
 [h_2\circ p_3](x)\cdot h_3(x)\Bigr)\,d\mu(x)
\nonumber
\end{gather}
We substitute $x=p_3(y)$ to (\ref{eq2.20})
and replace the summand (\ref{eq2.20}) by
$$
\int_X
c\Bigl([h_1\circ p_2\circ p_3](x)\,\,,\,\,
[h_2\circ p_3](x)\Bigr)\,d\mu(x)
$$
We obtain that (\ref{eq2.19}) equals (\ref{eq2.21}) due the
cocycle identity (\ref{eq2.3}) for $c(\cdot,\cdot)$ with
$$
g_1=h_1\circ p_2\circ p_3(x),\quad
g_2=h_2\circ p_3(x),
\quad
g_3=h_3(x)
$$

%

\smallskip

{\sc Remark.} To apply this construction,
we need in a group $G$, having a nontrivial
$\R$-central extension.
Some reasonable examples are
$G=\R/\Z$, $\U(p,q)$, $\mathrm {SO}^*(2n)$,
$\mathrm{SO}(n,2)$, and the group
$\mathrm{SDiff}(S^1)$ of orientation preserving diffeomorphisms of
the circle.

\smallskip


{\bf\punct \label{ss:another-fB} 
Another explanation of the central extension of $\fB(X,G)$.}
Denote by $\wt G$ the central extension of $G$
 defined by the cocycle $c$, by definition $\wt G\supset \R$. 
Concider the group
$\fB(X,\wt G)$ and its central Abelian subgroup 
$\cF(X,\R)$.    
Denote by $Q$ the subgroup of 
$\cF(X,\R)$ consisting of functions $f$ such that
$$
\int_X f(x)\,d\mu(x)=0
$$
Obviously, $Q$ is a normal subgroup in $\fB(X,\wt G)$.
It can readily checked that
$$
\wt\fB(X,G)=\fB(X,\wt G)/Q
$$


{\bf\punct Embedding $\Symp(M)\to\fB(M,\Sp(2n,\R))$.}
Now let $X\simeq M$ be a $2n$-dimensional
 symplectic manifold, let
$\Omega\subset\R^{2n}$ be an open domain,
and let $\iota:\Omega\to M$ be a symplectic
embedding such that the set $M\setminus\iota(\Omega)$
has a zero measure. For $g\in\Symp(M)$,
consider the map $q:=\iota^{-1}g\iota:\Omega\to\Omega$
defined almost sure.
Denote by $q'(x)$ the Jacobi matrix of $g$ at point $x$.

The map
$$g\mapsto (q,q')$$
is an embedding $\Symp(M)\to\fB(M,\Sp(2n,\R))$.
Our cocycle (\ref{eq1.5}) on $\Symp$
is induced from the cocycle (\ref{eq2.18})
on $\fB(M,\Sp(2n,\R)$.


{\bf \punct Groups of automorphisms
$\Aut(R,M)$ of bundles.\label{ss2.5}}
Now let $X=M$ be an  $m$-dimensional manifold
equipped with a volume form $\Omega$.
Consider an   $2n$-dimensional
real vector
bundle $R\to M$ on $M$
with fibers $R_x$, $x\in M$.
 Assume that fibers $R_x$
are equipped with skew-symmetric bilinear forms
$\{\cdot,\cdot\}_x$.

Denote by $\Aut(R,M)$ the group of smooth
maps $\Theta:R\to R$
 satisfying the conditions

 1. An image of any fiber $R_x$ is some fiber $R_{\theta (x)}$,
 and the map $\Theta$ induces a linear map $R_x$ to $R_{\theta (x)}$
 preserving the skew-symmetric form, i.e.,
 $$\{\Theta v, \Theta w\}_{\xi(x)}=\{v,w\}_{x}$$
 where $v$, $w\in R_x$.

  2. The map $x\mapsto \theta(x)$
  is a diffeomorphism of the base $M$
   preserving the volume form $\Omega$.

  3. For a noncompact manifold $M$, we have an additional conditions:
there is a compact subset $L\subset M$ such that for
any $x\notin L$ we have $\theta(x)=x$ and $\Theta v=v$
 for any $v\in\R_x$.

 \smallskip

The group $\Aut(R,M)$ admits an obvious embedding to the group
$\fB(M,\Sp(2n,\R))$. Indeed,  consider an arbitrary trivialization
of the bundle $R$ over an arbitrary open set $X\subset M$
 of a complete
measure. Then elements of $\Aut(R,M)$ induce transformations
of  $\cL(X,\R^{2n})$ of the type (\ref{eq2.17}).

The embedding $\Aut(R,M)\to \fB(M,\Sp(2n,\R))$ described above
is not canonical, but any two embeddings
$\eta_1$, $\eta_2$ of this type
are conjugated by some element $r\in\fB(M,\Sp(2n,\R))$:
\begin{equation}
\eta_2(\Theta)=r^{-1} \eta_1(\Theta) r
\label{eq2.23}
\end{equation}

Thus, the central extension of $\fB(M, \Sp(2n,\R)$ induces
a central extension of the group $\Aut(R,M)$.
The formula for the cocycle
 depends on the set $X\subset M$
 and on a trivialization of the bundle.
 But due (\ref{eq2.23}) all these cocycles are
 equivalent. Thus, our central extension is canonical.

\smallskip


{\bf \punct A geometric construction
of the central extension of $\Aut(R,M)$.\label{ss2.6}}
Let $R$, $M$ be the same as above.

Consider
an almost complex structure on the bundle $R$.
Recall that this is an operator $J_x:R_x\to R_x$
depending on $x$ smoothly and satisfying
$$
J_x^2=-1; \qquad \{J_x v,J_x w\}_x=\{v,w\}_x
$$
We also assume that the symmetric bilinear
 form
 $$S_x(v,w)=\{J_x v,w\}_x$$
 is positive definite on each fiber.
Such almost complex structures exist,
see, for instance, \cite{MS}
(evidently, such structure is not unique).
In particular, the tangent bundle $TM$ became a complex
$n$-dimensional bundle, but in this moment
we prefer to consider $TM$ is a real bundle with an additional
structure.

Consider the fiber-wise
 complexification $R^\C$ of the vector bundle
$R\to M$. In each fiber $(R_x)^\C$
of $R^\C$, we have two
canonically defined subspaces
$$R_x^\pm:=\ker (J_x\mp i),
\qquad (R_x)^\C=R_x^+\oplus R_x^-$$

Thus we can represent any real symplectic linear operator
$h:R_x\to R_y$ as a complex block operator
$$h: R_x^+\oplus R_x^-\to R_y^+ \oplus R_y^-$$
We denote by $\Phi(h)$ the block corresponding
$R_x^+\to R_y^+$.

The formula for the 2-cocycle is
$$
C(\Theta^{(1)},\Theta^{(2)})=
\int_M
c\Bigl((\Theta^{(1)}\Bigr|_{R_{\theta^{(2)}(x)}},
\Theta^{(2)}\Bigr|_{R_{x}}
\Bigr)\Omega(x)
$$


{\bf \punct Groups of symplectomorphisms.\label{ss2.7}}
Now let $M$ be a $2n$-dimensional
 symplectic manifold. Denote by $R$ the tangent
 bundle to $M$.
 We have the natural embedding
$$
\Symp(M)\subset \Aut(R,M)
$$
and hence we can induce a 2-cocycle
from the group $\Aut(R,M)$.

We also can put $\Symp(M)$ to $\fB(M,\Sp(2n,\R))$ directly, as
it was done above in \ref{ss1.3}.

The cohomological class of the induced cocycle
 do not depend on the embeddings
since all these embeddings are conjugated
by interior automorphisms
of $\fB(M,\Sp(2n,\R))$.

\smallskip


{\bf \punct Proof of Proposition \ref{prop1.2}.\label{ss2.8}}
Let $\Xi$ be a domain in $\R^{2n}$.
Then
\begin{equation}
\nu:(g,x)=\det\Phi(g'(x))/|\det\Phi(g'(x))|
\label{eq2.24}
\end{equation}
is a well defined function
$$
\SSymp(\Xi)\times \Xi \to \T
$$
where $\T$ is the group of complex numbers
$z$ such that $|z|=1$.

By the covering homotopy theorem, this map
can be lifted to the map
$$
\wt\nu: \SSymp^\sim(\Xi)\times \Xi\to \R
$$
such that
$$\exp\bigl(i \wt\nu(g,x)\bigr)=\nu(g,x),\qquad \wt\nu(e,x)=0$$
The function $\wt\nu$ is precisely a single-valued branch
of
$$
\Im\ln\det\Phi(g'(x))=\ln\Bigl[
      \det\Phi(g'(x))/|\det\Phi(g'(x))|\Bigr]
$$
Then
\begin{equation}
\Gamma^\circ(g):=\int_\Xi \wt\nu(g,x)\,dx
\label{eq2.25}
\end{equation}
is a trivializer of the 2-cocycle $C(\cdot,\cdot)$
on $\SSymp(\Xi)$, i.e.
$$
C(g_1,g_2)=\Gamma^\circ(g_1g_2)-\Gamma^\circ(g_1)-\Gamma^\circ(g_2)
$$
Indeed, the right hand side is
\begin{multline*}
\int \Im\ln \Phi\bigl([g_1'\circ g_2](x)\cdot g'_2(x)\bigr)\,dx
-\int \Im\ln \Phi(g'_1(x))\,dx
-\int \Im\ln\Phi (g'_2(x))\, dx
\end{multline*}
We change the variable $x\mapsto g_2(x)$ in the second summand
and transform the expression to the form
$$
\int\Im\tr\ln\Bigl[\Phi\bigl([g_1'\circ g_2](x) \bigr)^{-1}
          \Phi\bigl([g_1'\circ g_2](x)\cdot g_2'(x) \bigr)
          \Phi\bigl(g_2'(x)^{-1} \bigr]\Bigr]\,dx
$$

This proves Proposition \ref{prop1.2}.

{\sc Remark.}
Generally, this trivializer is not unique, since
(see \ref{ss2.1})
there are nontrivial homomorphisms
$\SSymp(\Xi)^\sim\to\R$, namely the flux homomorphisms
and the Calabi invariant (see, for instance \cite{Ban},
\cite{MS}, we discuss them below in a 
special case).
By the Banyaga Theorem \cite{Ban}, this list exhaust all
the  homomorphisms from $\SSymp^\sim$ to Abelian groups.

\smallskip

\begin{corollary}
\label{cor2.3}
 Let $\Xi\subset\R^4$ be an open domain.
Then our central extension  of the (disconnected)
group $\Symp(\Xi)$ is trivial.
\end{corollary}

\smallskip

Indeed, the group $\Symp(\R^4)$ is connected and contractible
(Gromov \cite{Gro}). By Proposition \ref{prop1.2},
 its central extension
is trivial. But $\Symp(\Xi)\subset \Symp(\R^4)$.

\smallskip


{\bf\punct \label{ss2.9} Proof of Proposition \ref{prop1.3}.}
In the case of general symplectic manifolds
$M$, the sense of the expression (\ref{eq2.24}) is not clear,
since an operator $\Phi(\cdot)$
maps one fiber of a bundle to another one
and hence its determinant is not well defined.

Consider an Hermitian structure on the tangent
bundle $T\,M$ to $M$ compatible with the symplectic form
(as in \ref{ss2.6}).
Then $T M$ becomes a complex bundle, denote it by
$\mathcal T M$ (it is also isomorphic to the subbundle
$R^+\subset R^\C$ defined above).
Consider
its maximal exterior power
$\wedge^n \mathcal T M$.

Assume that the  bundle
$\wedge^n \mathcal T M$
is trivial. Fix its trivialization.
This precisely means that we have defined
determinants of the operators $\Phi(\cdot)$ connecting
different fibers.

Then we repeat our argument with a covering homotopy
and obtain Proposition \ref{prop1.3}.

\bigskip

 \section{Two-dimensional surfaces}

\setcounter{equation}{0}

\bigskip

\nopagebreak

Here we prove nontriviality of the cocycle
$C$ in the case
of a two-dimensional oriented
surface $\M$ of genus $g\ge 3$.
Evidently, in 2-dimensional case
the group $\SSymp(\M)$ coincides with the group
of volume preserving and orientation preserving
diffeomorphisms.
Also, $\Sp(2,\R)=\SL(2,\R)$.

We fix the following notation

--- $\Lambda$ is a simply connected domain
in $\R^2$ equipped with the form
$$\omega=dx\wedge dy$$
 It convenient to think that $\Lambda$ is
a standard disk.

--- $\Delta$ is a multiconnected domain in $\R^2$.
It is convenient to think, that $\Delta$ is a disk
with $k>0$ holes and the exterior boundary of $\Delta$
is the boundary of $\Lambda$.

--- $\sigma(S)$ is the area of a domain $S\subset\R^2$.

\qquad

The group $\Symp(\Lambda)=\SSymp(\Lambda)$
 is connected and contractible
(Smale). The groups $\Symp(\Delta)$ are disconnected
(this is obvious).
The identity component $\SSymp(\Delta)$ is contractible
(Earle, Eells \cite{EE}).

\smallskip

Our main arguments for a proof are:
continuity of the cocycle in the topology of
convergence in measure, formula (\ref{eq2.25}) for
a trivializer in a flat case, a certain Dehn relation
in the Teichmuller group, and the Banyaga Theorem.

\smallskip


{\bf\punct \label{3.1} Global angle of rotation.}
For $\psi\in\R$ consider the unit vector
$$
v_\psi=\cos\psi e_1+\sin\psi e_2
$$
Denote by $S^1$ the set of all unit vectors
in $\R^2$.

Let $q\in\Symp(\Lambda)$.
Consider a point $x\in\Lambda$ and
a unit vector
$v_\psi$ applied at this point.
Consider the image $w=q'(x)v$ of $v$ under the Jacobi
matrix $q'(x)$. The normalized  vector
$w/\|w\|$ has a form $v_\phi$.
We assume
$$\ang(q,x,v):=\phi-\psi$$
i.e., $\ang(\cdot)$ is the angle of turning of a vector
under a diffeomorphism.

\smallskip

\begin{lemma}
 \label{l3.1}
 a)
 There exists a unique continuous function
(global turning angle)
$$\Ang:\SSymp(\Lambda)\times\Lambda\times S^1\to \R$$
such that

\smallskip

1)  The composition of $\Ang(\cdot)$ and the map
$\R\to\R/2\pi \Z$ is $\ang(\cdot)$.

2) $\Ang(e,x,v)=0$.

3)  Let $U$ be a neighborhood
of boundary of $\Lambda$, where $q(x)=x$.
 Then
  $$
  \Ang(q,x,v)=0\qquad \text {for all $x\in U$ and all $v\in S^1$}
  $$

 b)  The function $\Ang$ satisfies the identity
 \begin{equation}
 \Ang(q_1\circ q_2,x,v)=
\Ang(q_2,x,v)+\Ang(q_1,\,q_2(x),\,q_2'(x)v)
\label{eq3.1}
 \end{equation}
\end{lemma}

  {\sc Proof.}
a) Consider a map
$$f:\SSymp(\Lambda)\times\Lambda\times S^1\to S^1$$
defined in the following way:
$$f(q,x,v)=\frac{q'(x)v}{\|q'(x)v\|}$$
Next, we consider the covering map
$$
\wt f:\SSymp(\Lambda)\times\Lambda\times \R\to \R;
\qquad \wt f(e,x,\phi)=\phi
$$
defined by the covering homotopy Theorem
(recall that $\SSymp(\Lambda)$ is simply
connected).
Then
$$
\Ang(q,x,v_\psi)=\wt f(q,x,\psi)-\psi
$$

b) is obvious.
\hfill $\square$

\smallskip

\begin{corollary}
 \label{cor3.2} 
 The global turning angle
 is well defined on the
(disconnected) group
$\Symp(\Delta)$.
\end{corollary}

\smallskip

{\sc Proof.} Indeed, each
compactly supported symplectomorphism of
$\Delta\subset \Lambda$ is also a symplectomorphism of $\Lambda$.
\hfill $\square$

\smallskip

In particular, Lemma \ref{l3.1}
 gives the following geometrically
visual way of evaluation of $\Ang(\cdot)$.

\smallskip

\begin{lemma}
 \label{l3.3}
  Let $q\in\SSymp(\Lambda)$.
Consider a point $z$ near the boundary of $\Lambda$, where
$q(z)=z$. Let $\ell(t)$ be a smooth curve, $\ell(0)=z$,
$\ell(1)=x$, $\frac d{dt}\ell(1)=v$. Then
$$
\Ang(q,x,v)=\bigl\{\text{total turning of vector
$\frac d{dt} (q(\ell(t))$}\bigr\}-
\bigl\{\text{total turning of vector
$\frac d{dt} \ell(t)$}\bigr\}
$$
\end{lemma}

{\sc Proof.} Indeed, $\Ang(\cdot)$ must be continuous
along the curve
$$
\xi(t):=(q,\ell(t),\dot\ell(t))\in\SSymp\times\Lambda
\times S^1 \qquad\qquad\qquad\qquad\square
$$

Let $\Phi(h)$, where $h\in\SL(2,\R)$, be the same as above.
In our case $\Phi(h)$ is an element of $\C$
 and $|\Phi(h)|\ge 1$.
 Let
 \begin{equation}
 h=S A,\qquad
\label{eq3.2}
 \end{equation}
be the polar decomposition of $h$, where $A$ is a rotation
by some angle $\theta$
and $S$ is a contraction-dilatation with respect two orthogonal
 axes. Then
 \begin{equation}
 \Phi(h)/|\Phi(h)|=e^{i\theta}
\label{eq3.3}
 \end{equation}

\begin{lemma}
\label{l3.4}
  Consider a continuous branch
of the function
$$\gamma^\circ(g,x):=\Im\ln\Phi(g'(x))$$
on $\SSymp(\Lambda)\times\Lambda)$ such that
$\gamma(e,x)=0$.
We have
\begin{equation}
|\Ang(q,x,v)-\Im\ln\Phi(g'(x))|<\pi/2
\label{eq3.4}
\end{equation}
for all $q\in\SSymp(\Lambda)$,
$x\in\Lambda$, $v\in S^1$.
\end{lemma}


\smallskip

{\sc Proof.} We can define the both functions in the following
way. Consider the map
$$
F:\SSymp(\Lambda)\times\Lambda\to\SL(2,\R)
$$
given by
$(q,x)\mapsto q'(x)$.
Consider the covering map
$$\wt F:\SSymp(\Lambda)\times\Lambda\to\SL(2,\R)^\sim$$

Let $h$ ranges in $\SL(2,\R)$.
The function
$\alpha(h):=\Im\ln\Phi(h)$ is a function on $\SL(2,\R)^\sim$.
Also the $\R$-valued angle
of turning of a unit vector $v\in\R^2$
under $h\in\SL(2,\R)^\sim$
is a function on $\SL(2,\R)^\sim$
(denote it by $\beta_v(h)$).
 We have
 \begin{equation}
 \Ang(g,x,v)=\beta_v(\wt F(g'(x))),\qquad
 \gamma^\circ(g,x)=\alpha(\wt F(g'(x)))
\label{eq3.5}
 \end{equation}
 Hence it is sufficient to show that
\begin{equation}
 |\alpha(h)-\beta_v(h)|<\pi/2
\label{eq3.6}
 \end{equation}
 This is a corollary of the following
 statement:

 -- Let $h\in\SL(2,\R)$, let $h=S A $
 be its polar decomposition (\ref{eq3.2}).
  Then the angle between $hv$ and $Av$
 is less than $\pi/2$.

 The latter statement is obvious.
 Passing to the covering group $\SL(2,\R)^\sim$,
  we obtain (\ref{eq3.6});
 applying (\ref{eq3.5}), we obtain (\ref{eq3.4}).
\hfill$\square$

\smallskip

\begin{lemma}
 \label{l3.5}
$$
\qquad \max_{x,v} |\Ang(q_1\circ q_2,x,v)|
 \le  \max_{x,v} |\Ang(q_1,x,v)|+ \max_{x,v} |\Ang(q_2,x,v)|
$$
\end{lemma}
%

\smallskip

{\sc Proof.} This follows from (\ref{eq3.1}).

\smallskip

\begin{figure}
\begin{center}
\includegraphics[scale=0.5]{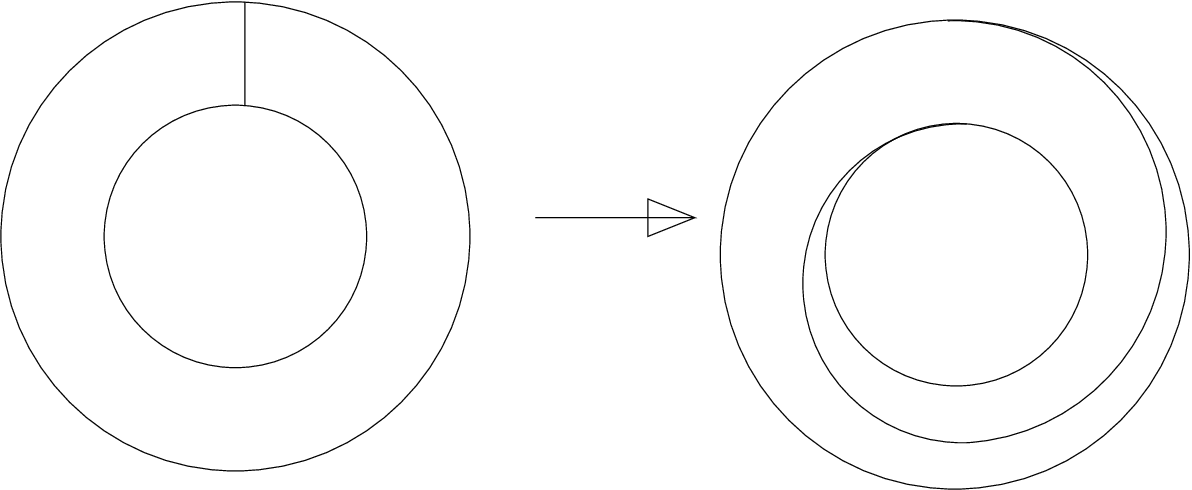}
\end{center}
\caption{ A standard twist.\label{fig1}}
\end{figure}

{\bf\punct \label{ss3.2} Twists.}
Let $r>0$, $\phi$ be the polar coordinates on the plane.
Consider a ring $S:a\le r\le b$.
A {\it standard right twist} is
a diffeomorphism $q:S\to S$
 that fix the both boundary circles,
 i.e.,
$q(a,\phi)=(a,\phi)$,
$q(b,\phi)=(b,\phi)$,
and
$$
q(r,\phi)=q(r,\phi+\mu(r))
$$
where $\mu$ is an arbitrary
 smooth increasing function
on $[0,\infty)$ such that
\begin{align*}
\mu(x)=0\qquad
\text{for $x<a+\delta$};\qquad
\mu(r)=2\pi\qquad \text{for $x>b-\delta$}
\end{align*}
for some $\delta>0$, see Fig. \ref{fig1}.

{\sc Remark.} The image of a right
twist under the orientation
preserving map $(r,\phi)\mapsto (ab/r,-\phi)$
is a right twist again. A diffeomorphism inverse to
a right twist is  not a right twist (it is called
a {\it left twist}).
\hfill$\square$

Consider a closed smooth non-self-intersecting
 curve $C$ on the surface $\M$. Consider a 'small' neighborhood
 $U$
 of $C$. For some standard
 ring $S\subset \R^2$
 consider    some area-preserving
 and orientation preserving diffeomorphism
  $p:S\to U$.
Consider a diffeomorphism $h\in\Symp(\M)$
having the form

-- $h(m)=m$ if $m\notin U$

-- $p^{-1} h p$ is a standard right twist of $S$.

We call such diffeomorphism by a {\it twist about the curve} $C$
supported by the neighborhood $U$.

A {\it Dehn twist} is a twist about a nonseparating curve
$B$ (i.e., the domain $\M\setminus B$ is connected).
For any two nonseparating curves $B_1$, $B_2$
there is a diffeomorphism $q\in\Symp(\M)$
such that $q(B_1)=B_2$ (see, for instance,
\cite{Iva}). If $h$ is a right twist about $B_1$,
then $q^{-1} h q$ is a right twist about $B_2$.

\smallskip


{\bf\punct \label{ss3.3} $\epsilon$-twists.}
Below we use families of twists depending on parameters
and we are need in some uniform estimates in parameters.
By this reason, we give more rigid definitions.
%

\SS

Fix a function $\nu$ on $(-\infty,\infty)$
satisfying the conditions

\SS

-- $\nu=0$ on the ray $(-\infty,0)$
and $\nu=2\pi$ on the ray $x\ge 1/2$

\SS

-- $\nu$ is $C^\infty$-smooth and increasing.

\SS

 This function remains fixed until the end of this section.

\SS

A {\it standard $\epsilon$-twist}
is a family $q(\epsilon)$ of diffeomorphisms of the ring
$1/2<r<3/2$ depending in a parameter $\epsilon$,
they  have the form
$$
\bigl(r,\phi\bigr)
\mapsto \bigl(r, \phi+\nu\bigl(\tfrac{r-1}\epsilon\bigr)\bigr),
 \qquad 0<\epsilon<1/2
$$
(in particular, this map is identical outside the ring
$1<r<1+\epsilon/2$).

\SS

Now let $C$ be a closed non-self-intersecting
 curve on $\M$. Let $h$ be an embedding
 of some ring $1-\delta<r<1+\delta$
 to $\M$ such that the $C$ is the image of the circle
$r=1$.
 We define an
  {\it $\epsilon$-twist} about $C$
  as a family of diffeomorphisms of $\M$ depending on
  $\epsilon$ and given by
\begin{align*}
m&\mapsto h\circ q(\epsilon)\circ h^{-1} (m),
   \qquad \text{ for $m\in U$}
\\
m&\mapsto m,\qquad \text{ for $m\notin U$}
\end{align*}
 The parameter $\epsilon$ ranges
  in the interval $(0,\delta)$.

\SS

 We fix some notation

\SS

  -- $\fT_C(\epsilon)$ for an $\epsilon$-twist about $C$;
 we omit $h$ from our notation, but we remember that $h$
  is fixed.

\SS

  -- $\fU(\epsilon)=\fU_C(\epsilon)$
for the support of the twist $\fT_C(\epsilon)$.
The area of $\fU_C(\epsilon)$
 tends to 0 as $\epsilon\to 0$, more precisely,
 $\sigma(\fU_C(\epsilon))=O(\epsilon)$.

\smallskip


{\bf\punct \label{ss3.4} 
Trivializator of the cocycle $C$ on a twist.}
The cocycle $C(\cdot,\cdot)$ on the group $\SSymp(\Lambda)$
is trivial,  its canonical trivializer
\begin{equation}
\Gamma^\circ(q)=\int_{\Lambda}
\Im\ln\Phi(q'(x))\,dx
\label{eq3.7}
\end{equation}
was defined above in \ref{ss2.8}.

\smallskip

\begin{lemma}
\label{l3.6}
 Let $B$ be a smooth non
 self-intersecting curve in $\Lambda$
surrounding the domain $S$. Then
$$
\Gamma^\circ(\fT_B(\epsilon))=2\pi \sigma(S)+O(\epsilon),
\qquad \text{as $\epsilon\to 0$}
$$
\end{lemma}

{\sc Proof.} The strip $\fU(\epsilon)$ separates $\Lambda$
into two domains, i.e., the exterior domain $W^{ext}$
and the interior domain $W^{int}$.
Let $\ell(t)$ be a (short) curve intersecting
the strip $\fU(\epsilon)$, let $\ell(0)\in W^{ext}$,
$\ell(\delta)\in W^{int}$. The Jacobi matrix
$[\fT_B'\circ\ell](t)$ is 1
at $t=0$ and at $t=\delta$. But the curve
$[\fT'\circ\ell](t)$ is noncontractible in $\SL(2,(\R))$
and moreover, it is
 a generator of the fundamental group $\pi_1(\SL_2(\R))$
 ( it is
 sufficient to verify that the vector
 $[\fT_B'\circ\ell](t)\dot\ell(t)$ is turning
  by the angle $2\pi$ as we pass $\ell(t)$).

 Hence $\ln\Phi(\fT'_B(\ell(\delta))=2\pi$,
 and thus $\ln\Phi(\fT'_B(x))=2\pi$ for $x\in W^{int}$.

 Further,
$$
\Gamma^\circ(\fT_B(\epsilon))=\int_{\Lambda}
\Im\ln\Phi(\fT_B'(x))\,dx=\int_{W^{ext}}+\int_{\fU(\epsilon)}
+\int_{W^{int}}
$$
The integrand in the first summand is 0,
the integrand in the last summand is $2\pi$
(and hence the integral over $W^{int}$ is
$2\pi\sigma(S)+O(\epsilon)$).
Thus it is sufficient to show that
$\gamma^\circ=\Im\ln\Phi(\fT'_B(x))$ is bounded in
the thin strip $\fU(\epsilon)$.
By Lemma \ref{l3.4}, it is sufficient
to show the boundedness of $\Ang(\fT_B, x, v)$,
and by Lemma \ref{l3.3}, the latter statement is obvious.
\hfill $\square$.


\smallskip

{ \bf\punct \label{ss3.5} Flux homomorphisms.}
Now, let $\Lambda\subset\R^2$ be a disk (or simply connected domain)
bounded by a curve $Z_0$.
Let  $\Delta\subset\Lambda$ be a multi-connected domain
bounded by smooth curves $Z_0$ (the exterior boundary)
and
$Z_1$, \dots, $Z_k$
(boundaries of the holes).
Denote by $\omega$ the standard symplectic form
$dx\wedge dy$ on $\R^2$. Let $\lambda$ be a $1$-form
on $\R^2$ such that $d\lambda=\omega$
(for instanse, $\lambda= x\,dy$).


For each $j$,
fix a non self-intersecting
 curve $u_j(t)$ connecting $Z_0$ and $Z_j$,
 let $u_j(0)\in Z_0$, $u_j(1)\in Z_j$.
 We define the function $\tau_j(g)$ in
the variable $g\in\Symp(\Delta)$ by 
\begin{equation}
\tau_j(g)=\int_{u_j} \lambda-\int_{u_j} g^*\lambda
=\int_{u_j}\lambda-\int_{gu_j}\lambda=
\int_{D_j}\omega
\label{eq:flux}
\end{equation}
where $D_j$ is a 2-cycle whose boundary is
$u_j-gu_j$.
Non-formally,
 $\tau_j(g)$ is the oriented area
 of the domain bounded by the curves
$u_j(t)$
 and $g(u_j(t))$.

\begin{lemma}
\label{l:flux}
a) $\tau_j$ does not depend on a choice of $u_j$.

b) $\tau_j$ is a homomorphism $\Symp(\Delta)\to\R$. 
\end{lemma}

   The maps $\tau_j$ are called by 
{\it flux homomorphisms}, for definitions
and properties in a general symplectic case,
see Banyaga \cite{Ban}, McDuff, Salamon \cite{MS}.

\SS

{\sc Proof.} a) Let $u_j'$ be another curve, let $\tau_j'$ be another 
map. Let $R$ be a 2-cycle on $\R^2$ with boundary 
$u_j-u_j'$,
$$
\tau_j(g)-\tau_j'(g)=
\Bigl(\int_{u_j}\lambda-\int_{u_j'}\lambda)\Bigr)-
\Bigl(\int_{gu_j}\lambda-\int_{gu_j'}\lambda\Bigr) 
=\int_R \omega -\int_{gR}\omega=0
$$

b) 
$$
\tau_j(g_1 g_2)=\int_{u_j}\lambda- \int_{g_1g_2u_j}\lambda
=
\int_{u_j}\lambda-\int_{g_2u}\lambda+ \int_{g_2u}\lambda-
 \int_{g_1(g_2u_j)}\lambda=
\tau_j(g_2)+\tau_j(g_1)
$$





\smallskip

{\sc Remark.} In a general case, flux homomorphisms
are defined on the connected group $\SSymp$. In our case,
we obtain homomorphisms $\tau_j$ of group $\Symp(\Delta)$,
these homomorphisms depend on an embedding of the symplectic
manifold $\Delta$ to $\R^2$, since the areas of the holes
$Z_j$ participate in formula (\ref{eq:flux}) 


\SS

{\bf \punct Values of fluxes on twists.}
We preserve the notation of
the previous subsection. The following
statement is obvious. 

\begin{lemma}
\label{l:twist-flux}
Let $C$ be a Jordan contour in $\Delta$ surrounding
a domain $S$. Then 
$$
\tau_j(\fT_C(\epsilon))=
\begin{cases}
 \sigma(S)+O(\epsilon),\,\,
  &\text{if $Z_j\subset S$}
\\
0,&\text{otherwise}
\end{cases}
$$
\end{lemma}


{\bf\punct \label{ss3.6} The Calabi homomorphism.}
We preserve the notation of the previous subsection.
Consider a 1-form $\lambda$ on $\R^2$ such that
$d\lambda=\omega$.

For $g\in \SSymp(\Lambda)$,
the form $g^*\lambda-\lambda$ is
closed and  hence it is exact.

Hence
\begin{equation}
g^*\lambda-\lambda=d F
\label{eq3.8}
\end{equation}
 The function $F$ is defined up to
an additive constant. We assume $F=0$ on $Z_0$.
Then
\begin{equation}
\kappa(g):= \int_{\Lambda} F dx\,dy
\label{eq3.9}
\end{equation}
is a homomorphism $\SSymp(\Lambda)\to \R$.
It can easily be checked that $\kappa(g)$
does not depend on a choice of a potential $\lambda$.
The homomorphism $\kappa(g)$
 is called by the {\it Calabi invariant}.

\SS

 Next, we  restrict the Calabi invariant to
 the group
 $\Symp(\Delta)\subset\SSymp(\Lambda)$
 and hence we obtain the homomorphisms
 $
 \Symp(\Delta)\to\R
 $.

 {\sc Remark.}
 In general situation, the
Calabi invariant is defined on  the kernel of all the flux
homomorphisms $\subset\Symp$.
 In our case, it is defined globally,
but its definition is not canonical, it depends on
an embedding of the symplectic domain $\Delta$ to $\R^2$.
Indeed, formula (\ref{eq3.9}) includes
  integration over the holes,
this operation is not invariantly defined on the manifold
$\Delta$.
But these integrals over the holes give 
a linear combination of flux homomorphisms.

\smallskip

By the Banyaga Theorem (see \cite{Ban}, see also
some useful additions in Rousseau \cite{Rou},
see also \cite{Ban2}),
the intersection of the kernels of all the flux
homomorphisms and of the kernel
of the Calabi invariant is a perfect group
(i.e., it has no homomorphisms to Abelian groups;
moreover, this intersection is a simple group).

Thus, in our case, each measurable%
\footnote{The Choice Axiom implies an "existence"
of  non-measurable homomoprhisms
$\psi:\R\to\R$ (number of such homomorphisms is
$2^{continuum}$). A composition of $\psi$ and a flux homomorphism
is non-measurable homomorphism $\SSymp\to\R$.}
 homomorphism $\SSymp(\Delta)\to \R$
is a linear combination
of $k$ flux homomorphisms
$\tau_j$ and of the Calabi
invariant $\kappa$.

\smallskip


{\bf \punct\label{ss3.7}
 Values of the Calabi invariants on twists.}

\nopagebreak

\smallskip

\begin{lemma}
\label{l3.7}
 Let $B$ be a smooth non
 self-intersecting curve in $\Lambda$
surrounding the domain $S$. Then
\begin{equation}
\kappa(\fT_B(\epsilon))=\sigma(S)^2+O(\epsilon),
\qquad\qquad \epsilon\to 0
\label{eq3.10}
\end{equation}
\end{lemma}

\smallskip

{\sc Proof.} We preserve the notation from the proof
of  previous Lemma \ref{l3.6}.
Let us choose
$\lambda=\frac 12 (x\,dy-y\,dx)$.
Let $F$ be the  function (\ref{eq3.8}).
  $$
  F(\ell(1))=
  \int_\ell (\fT_B^*\lambda-\lambda)=
  \int_{\fT_B\ell}\lambda-\int_\ell \lambda
  $$
We obtain an integral over a closed curve
$Q$
composed from $\ell$ and $\fT_B\ell$.

This curve lies in the strip $\fU_B(\epsilon)$
and surrounds $W^{int}$.
By the  Green formula,
 $F\simeq\sigma(S)$ in $W^{int}$.

We have
$$
\kappa(\fT_B(\epsilon))=
\int_\Lambda F\,dx\,dy=\int_{W^{int}}+\int_{W^{ext}}
+\int_{\fU_B(\epsilon)}
$$
The fist summand gives $\sigma(S)^2+O(\epsilon)$,
the second summand is 0. It remains
to show that $F$ is uniformly bounded in $x$ and $\epsilon$
in the strip $\fU(\epsilon)$.


The value of $F(\ell(s))$ for $0<s<1$
is the integral of $\lambda$ over a nonclosed
curve $L$ composed from $\ell(t)$ and $\fT(\ell(t))$
with $0<t<s$ ($\ell(t)$ is passed in the inverse direction).
We include this curve into a closed contour $C$
adding the direct segments $[0,\ell(0)]$,
$[0,\ell(s)]$. The integral of $\lambda$ over these segments
vanishes. Hence, by the Green formula, $F(\ell(s))$
is
 the oriented
area of the curvilinear  sector bounded by
the contour $C$.
Obviously, this area is uniformly bounded.
\hfill$\square$

\smallskip


{\bf\punct \label{ss3.8} Preliminary remarks on trivializers.} 
First, let us consider a domain $\Omega\subset\R^2$
and a map $\iota:\Omega\to\M$
as in \ref{ss1.3}. This allows to fix an explicit expression
(\ref{eq1.5}) for
the cocycle $C$. Below we will choose
$\Omega$ and $\iota$ in a certain appropriate way.

Assume, that our cocycle $C(q_1,q_2)$ is trivial
on $\Symp(\M)=\Symp(\M,\Omega,\iota)$,
 let $\Gamma(q)$ be its trivializer.
In other words, consider the space
$\Symp(\M,\Omega,\iota)\times\R$ with
the multiplication (\ref{eq1.6}).
For each $q\in\Symp(\M,\Omega,\iota)$, consider
the element
$$\widetilde q:=(q,\Gamma(q))\in
 \Symp(\M,\Omega,\iota)\times\R $$
Then
$$
\wt q_1 \wt q_2=\wt{q_1 q_2}
$$

\smallskip



\begin{lemma}
\label{l3.8} 
 For a diffeomorphism
$q\in \Symp(\M,\Omega,\iota)$ consider
the set $\Move(q)$
of all $x\in\Omega$ such that $q(x)\ne x$.
Then, for each $r\in \Symp(\M,\Omega,\iota)$,
$$
|C(q,r)|<\frac\pi2 \sigma(\Move(q)),
\qquad|C(r,q)|<\frac \pi 2 \sigma(\Move(q))
$$
\end{lemma}

In particular, the value of the cocycle
$C(q_1,q_2)$ is $O(\epsilon)$
if one of the arguments $q_1$, $q_2$ is an
$\epsilon$-twist.

\smallskip

{\sc Proof.} Let $g_1,g_2\in\Sp(2n,\R)$.
If $g_1=1$ or $g_2=1$, then $c(g_1,g_2)=0$,
see formula (\ref{eq2.10}). It remains to apply
Theorem \ref{th2.1}.d.

\smallskip

\begin{corollary}
\label{cor:trivializers}
$$
\bigl|\Gamma(g_1 \dots g_k)-(\Gamma(g_1)+\dots+\Gamma(g_k))\bigr|
\le \frac\pi 2 \sum \sigma_i(\Move(g_i))
$$
\end{corollary}



\smallskip

\begin{lemma}
\label{l3.9} 
There is a constant $H$ such that
$$\Gamma(\fT_C(\epsilon))=H+O(\epsilon)
,\qquad \epsilon \to 0$$
for each Dehn $\epsilon$-twist.
\end{lemma}

\smallskip

{\sc Proof.} {\it Asymptotics for a single twist.}
Consider a nonseparating non self-intersecting
curve $C$ on the surface. Consider a ring 
$\Theta\supset C$, identify this ring with a subdomain
in the flat circle $\Lambda$.
Our cocycle must be trivial on the group $\Symp(\Theta)$
 of symplectomorphisms of $\Theta$.    
On another side, our central extension admits a canonical
trivialization $\Gamma^\circ$ on the group 
$\SSymp(\Lambda)\supset \Symp(\Theta)$. Hence
$$
\Gamma(q)
=
\Gamma^\circ(q)+ u(q),\qquad q\in \Symp(\Theta)
$$
where $u$ is a homomorphism $\Symp(\Theta)\to\R)$.

The group $\Symp(\Theta)$ is a semidirect product
$\Z\ltimes\SSymp(\Theta)$. Let $\xi\in\Symp(\Theta)$
be a generator of the mapping class group 
$\Symp(\Theta)/\SSymp(\Theta)$. Then $q=\xi^n\circ q^*$,
where $q^*\in\SSymp(\Theta)$. 
Hence the homomorphism $u(q)$ must have a form
$$
u(q)=n\cdot u(\xi)+ a\tau(q^*)+b\kappa(q)
$$
Now let $q=\fT_C(\epsilon)$. Then all the terms
of $\Gamma(q)$ have asymptotics of the form 
$\mathrm{const}+O(\epsilon)$ (for $\Gamma^\circ$, see
Lemma \ref{l3.6}; for $\kappa$, see Lemma \ref{l3.7};
 for $\tau$
this is more-or-less obvious).

\smallskip

{\it Coincidence of asymtotics for different twists}.
Let $B$, $C$ be two nonseparating
non self-intersecting curves on $\M$.
Let us show, that there exists $ r\in\Symp(\M)$
such that

$$
\fT_C(\epsilon)= r\fT_B(\epsilon) r^{-1}
$$
for sufficiently small $\epsilon$.
Indeed, by definition of  $\epsilon$-twists,
we have fixed diffeomorphisms $h_B$, $h_C$
from some ring $1-\delta<r<1+\delta$
to some strips $\fU_B(\delta)$, $\fU_C(\delta)$
about $B$, $C$.
  For $m\in\fU_B(\delta)$,
we define $r$ as $r(m)=h_C\circ h_B^{-1}(m)$.
After this we extend $r$ to $\M\setminus \fU_B$
in an arbitrary way.
This is possible since $\M\setminus\fU_B$
and $\M\setminus\fU_B$ are symplectomorphic.

Now consider the corresponding elements
of $\Symp(\M,\Omega,\iota)$.
We preserve for them the same notation.
We have
$$
\wt \fT_C(\epsilon)\, \wt r=
\wt {\fT_C(\epsilon)  r}=\wt{r\fT_B(\epsilon)}= \wt r\,\wt\fT_B(\epsilon)
$$
Evaluating the trivializer
for the first and the last terms of this chain
we obtain
$$
\Gamma( \fT_C(\epsilon))+\Gamma( r)+
C(\fT_C(\epsilon),r)= \Gamma(r) +\Gamma(\fT_B(\epsilon))
+C(r,\fT_B(\epsilon))
$$
By Lemma \ref{l3.8}, we have
$\Gamma( \fT_C(\epsilon))-\Gamma( \fT_B(\epsilon))=O(\epsilon)$.
\hfill$\square$

\smallskip

{\bf\punct Proof of Theorem \ref{th1.4}.} 
Now consider an open domain $\widehat\Delta$ on $\M$
homeomorphic to a disk with 3 holes. Assume
also that the set $\M\setminus \widehat\Delta$ is connected.
Since the genus is $\ge 3$, this is possible.
Now we define more precisely the set $\Omega\subset\R^2$
and the map $\iota$. Let $\iota$ identify
symplectomorphically $\widehat\Delta$
and some disk $\Delta$ with 3 holes on $\R^2$;
denote by $\Lambda$ the simply connected domain
inside the exterior boundary of $\Delta$.

 On
$\M\setminus\widehat\Delta$ we can choose the map $\iota$
in an arbitrary way.

\begin{figure}
\begin{center}
\includegraphics[scale=0.5]{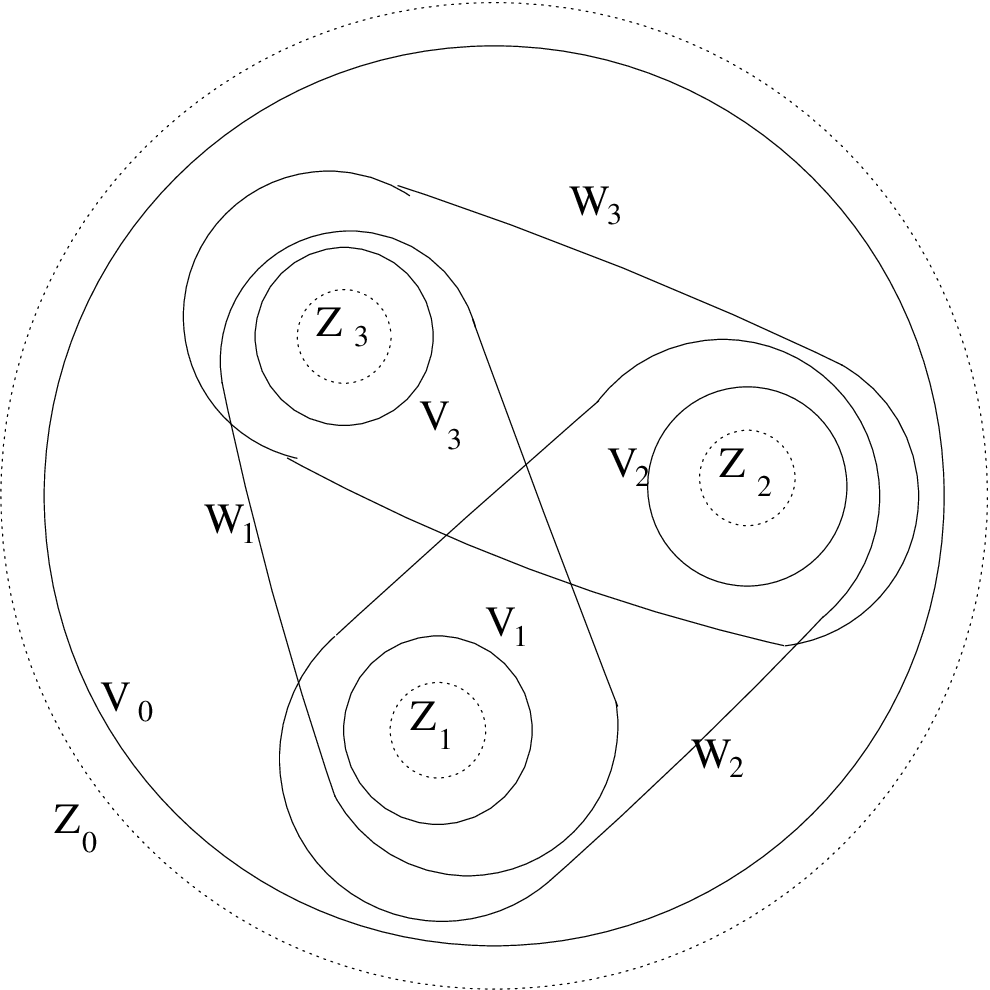}
\end{center}
\caption{Notation for curves.
The domain $\Delta$ is bounded 4 dotted
circles $Z_j$. The domain $\Lambda$ is the disk
bounded by $Z_0$.\label{fig3}}
\end{figure}

Consider 7 curves $V_0$, $V_1$, $V_2$, $V_3$,
$W_1$, $W_2$, $W_3$
 as on Fig.\ref{fig3}.

Consider the corresponding $\epsilon$-twists.
The diffeomorphism
$$
p(\epsilon):=\fT_{V_0}\fT_{V_1}\fT_{V_2}\fT_{V_3}
\fT_{W_1}^{-1}\fT_{W_2}^{-1}
\fT_{W_3}^{-1}
$$
is isotopic to the identity map, and moreover,
the isotopy can be done inside the domain $\Delta$;
this statement is the Dehn's {\it latern relation}
in Teichmuller group
rediscovered by Johnson,
 see \cite{Deh}, \cite{Joh},
\cite{Iva}).\hfill$\square$

\SS

All our curves are nonseparating in $\M$, and
by Lemma \ref{l3.9}
the trivializer $\Gamma$ on each of our 7 twists is
$H+O(\epsilon)$.
By Lemma \ref{l3.8},
$$
\Gamma(p)=H+O(\epsilon),\qquad \epsilon\to 0
$$
In particular, this means that
the main tern of the asymptotics of
$\Gamma(p(\epsilon))$ is invariant with respect to deformations
of our seven curves $V_0$,\dots.

\smallskip

On another hand, we have the trivializer
$\Gamma^\circ$
of $C(\cdot,\cdot)$ on $\SSymp(\Delta)$ described
in \ref{ss2.8} and given by (\ref{eq3.7}).
The difference of  two trivializers is a homomorphism
$\SSymp(\Delta)\to\R$.
 Hence it must
be a linear combination of 3 flux homomorphisms
$\tau_j$ and the Calabi
invariant.
Thus, we have
\begin{equation}
-H+
\Gamma^\circ(p)+\sum_{j=1}^3 a_j \tau_j(p)
+b\cdot\kappa(p)=O(\epsilon)
\label{eq3.11}
\end{equation}
where $a_j$, $b$ are some real constants.
We evaluate $\tau_j(p)$ and $\kappa(p)$ 
using Lemmae \ref{l:twist-flux}
and \ref{l3.7}. Denote by $\sigma[V_i]$ 
(resp. $\sigma[W_j]$) the area
surrounded by a contour $V_i$ (resp. $W_j$).
Then
\begin{multline*}
-H+ 
2\pi\sum_i \sigma[V_i]-2\pi\sum_j\sigma[W_j]+
\\
+
a_1(\sigma[V_0]+\sigma[V_1]-\sigma[W_1]-\sigma[W_2])
+
a_2(\sigma[V_0]+\sigma[V_2]-\sigma[W_2]-\sigma[W_3])+\\+
a_3(\sigma[V_0]+\sigma[V_3]-\sigma[W_3]-\sigma[W_1])
+\\+
b\bigl(\sum_i \sigma[V_i]^2-\sum_j\sigma[W_j]^2\bigr)
=O(\epsilon)
\end{multline*}
We can vary the areas $\sigma[\cdot]$
in arbitrary way in certain small intervals.
This is a contradiction.

\section{Nontriviality of the extension in the case of tori}

\setcounter{equation}{0}


{\bf\punct \label{ss4.1} Formulation of result.}
Consider the torus $\T^{2n}=\R^{2n}/\Z^{2n}$.
In notation of \ref{ss1.3},
it is natural to consider $\Omega=(0,1)^{2n}$
and the identical embedding $\Omega\to\R^n\to\T^{2n}$.
Then Theorem \ref{th1.5}
is a corollary of the following theorem.

\smallskip

\begin{theorem} 
\label{th4.1}
   For $n>2$, the $\R$-valued cocycle $c\in H^2(\Sp(2n,\R))$
is nontrivial on the group $\Sp(2n,\Z)$.
\end{theorem}

 Fix $\alpha>0$. Consider the homomorphism
 $\R\to \R/2\pi\alpha \Z$ and consider the image
 $c_\alpha$ of $c$ under this map.

 \smallskip

\begin{theorem}
\label{th4.2}  For $n>2$
and a noninteger $\alpha>0$, the
$\R/2\pi \alpha\Z$-valued cocycle $c_\alpha$
is nontrivial on the group $\Sp(2n,\Z)$.
\end{theorem}

\smallskip

Theorem \ref{th4.2} implies Theorem \ref{th4.1},
thus it is sufficient to prove Theorem \ref{th4.2}.

\smallskip

{\sc Remark.} For integer $\alpha$ the cocycle
$c_\alpha$ is trivial on the whole group $\Sp(2n,\R)$.

\smallskip

{\sc Remark.} Consider the subgroup 
$\Gamma_{1,2}\subset\Sp(2n,\Z)$
consisting of matrices 
$g=\begin{pmatrix}A&B\\C&D\end{pmatrix}$
such that diagonals entries of the matrices
$A^tC$, $B^tD$ are even. The cocycle
$c_{1/2}$
is trivial on $\Gamma_{1,2}$ (see \cite{Mum}, Section II.5),

\smallskip


{\bf\punct \label{ss4.2} Realization of the cocycle
$c_\alpha$ in a representation.}
Denote by $\B_n$ the set of
complex symmetric
($z=z^t$) matrices with norm $< 1$
(the Cartan matrix ball).

The symplectic group
acts on $\B_n$ by the linear-fractional transformations
$$z\mapsto z^{[g]}=(\Phi+z\ov\Psi)^{-1}(\Psi+ z\ov\Phi)$$

Fix $\alpha>0$.
Consider the representation $\widetilde T_\alpha$ of $\Sp(2n,\R)$
in the space of holomorphic functions on $\B_n$
given by
\begin{equation}
\widetilde
T_\alpha(g)f(z)=
  f(z^{[g]})\det(\Phi+ z\ov\Psi)^{-\alpha}
\label{eq4.1}
\end{equation}
It is a standard formula for a highest
weight representation of the group
$\Sp(2n,\R)$.

For $z\in \B_n$, the matrix
$\Phi+ z\ov\Psi$ is nondegenerat, see (\ref{eq2.14}).
Hence the expression
$$\det(\Phi+ z\ov\Psi)^{-\alpha}=
\det \Phi^{-\alpha} \det (1+z\ov \Psi \Phi^{-1})^{-\alpha}$$
has countable number of branches;
they are enumerated by values of $(\det\Phi)^{-\alpha}$.
Thus $\widetilde T_\alpha$ is a projective representation
of $\Sp(2n,\R)$
(or representation of the universal covering group $\Sp(2n,\R)^\sim$).

%



Nevertheless, we interpret the standard formula
(\ref{eq4.1}) in the following slightly
 nonstandard way.
Let us consider the normalized operators
$T_\alpha(g)$ given by
$$
T_\alpha (g)=
f(z^{[g]})\det(1+ z\ov\Psi\Phi^{-1})^{-\alpha}
$$
The expression
$$(1+ z\ov\Psi\Phi^{-1})^{-\alpha}=
\sum\limits_{k=0}^\infty \frac{(\alpha)_k z^k}{k!}
\bigl(-\ov\Psi\Phi^{-1}\bigr)^k
$$
is well defined as a sum of series.
Thus its determinant is well defined.

\smallskip

\begin{proposition}
 \label{prop4.3}.
\begin{equation}
T_\alpha(g_1)T_\alpha(g_2)=\sigma_\alpha(g_1, g_2)
     T_\alpha(g_1g_2)
\label{eq4.2}
\end{equation}
 where
\begin{equation}
\sigma_\alpha(g_1,g_2)=
 \det\Bigl[\Phi(g_1)^{-1}\Phi(g_1g_2)\Phi(g_2)^{-1}\Bigr]^{-\alpha}
= \exp\{-\alpha c(g_1,g_2)\}
\label{eq4.3}
\end{equation}
\end{proposition}

{\sc Proof.}
For $g_1$, $g_2$ near 1, it is proved by a trivial calculation.
After this we consider the analytic continuation
in $g_1$, $g_2$.
\hfill $\square$

The cocycles $\sigma_\alpha$ are precisely the
cocycles $c_\alpha$ in multiplicative notation.
We intend to prove Theorem \ref{th4.2} in the following form:

{\it the  restriction of the representation $T_\alpha$
to $\Sp(2n,\Z)$
can not be reduced to a linear representation
of $\Sp(2n,\Z)$
by a correction of the form
$$
T_\alpha(g)\mapsto \gamma(g) T_\alpha,\qquad \gamma(g)\in \C^*
$$
where $g$ ranges in $\Sp(2n,\Z)$.}

\smallskip


{\bf\punct \label{ss4.3} Another model of the same representation.}
By $\W_n$ we denote the Siegel wedge,
i.e., the set of complex symmetric matrices satisfying
$z$ with $\Im z>0$.

The group $\Sp(2n,\R)$ acts on $W_n$ by the transformations
$$z\mapsto z^{[g]}:= (a+zc)^{-1}(b+zd)$$
where
$g=\begin{pmatrix}
a&b\\c&d\end{pmatrix}
$
is a real symplectic matrix in the usual notation.
The action of $\Sp(2n,\R)$
 on $\W_n$ is given by
\begin{equation}
S_\alpha(g)f(z)=f(z^{[g]})\det(a+zc)^{-\alpha}
\label{eq4.4}
\end{equation}
It is well-known, that the (projective)
representation $S_\alpha$ is equivalent
to representation $T_\alpha$ defined above.
The intertwining operator
 is given by the transformation
$$
T_\alpha\Bigl[\frac 1{\sqrt 2}
\begin{pmatrix}
1&i\\i&1
\end{pmatrix}
\Bigr]f(z):=f((1+iz)^{-1}(i+z))
\cdot\det(1+iz)^{-\alpha}\cdot 2^{\alpha/2}
$$
where $z\in B_n$ (then its Cayley transform
 $(1+iz)^{-1}(i+z)$ is an element of $W_n$).
 

\SS

{\bf \punct\label{ss4.4} Linearization of representation
on a upper triangular subgroup.}
Consider the subgroup $\BB(\Z)\subset \Sp(2n,\Z)$ consisting
of the matrices
$$
\begin{pmatrix}
A&B\\0&A^{t-1}
\end{pmatrix}      ;\qquad\det A=1
$$

The cocycle $\sigma_\alpha$  is equivalent to trivial cocycle
on this subgroup. Indeed, the operators
\begin{equation}
S_\alpha
\begin{pmatrix}
A&B\\0&A^{t-1}
\end{pmatrix}
f(z)=f\bigl(A^{-1}(B+zA^{t-1})\bigr)
\label{eq4.5}
\end{equation}
define a linear representation of $\BB(\Z)$.

\smallskip

\begin{lemma} 
\label{l4.4} 
 Formula (\ref{eq4.5}) gives a unique
possible linearization of the representation
$S_\alpha$ on the subgroup $\BB(\Z)$.
\end{lemma}

\smallskip

This follows from the next lemma.

\smallskip

\begin{lemma}
 \label{l4.5} 
 The group $\BB(\Z)$
 has no homomorphisms to $\C^*$.
In particular, it has no homomorphisms to $\Z$ and $\Z_k$.
\end{lemma}

\smallskip

{\sc Proof.} Let $\chi:\BB(\Z)\to\C^*$ be a homomorphism.

First, the group $\BB(\Z)$ contains
the group $\SL(n,\Z)$
consisting of matrices
$\begin{pmatrix}A&0\\0&A^{t-1}\end{pmatrix}$.
This
group  has no Abelian quotients.
Hence, $\chi=1$ on $\SL(n,\Z)$.

 Also, $\BB(\Z)$ contains the group $N(\Z)$ consisting
of matrices
$\nu(B)= \begin{pmatrix} 1&B\\0&1
\end{pmatrix}$;
the product in this group corresponds to
the sum of the matrices $B$:
$$ 
\begin{pmatrix} 1&B_1\\0&1
\end{pmatrix}
\begin{pmatrix} 1&B_2\\0&1
\end{pmatrix}
=\begin{pmatrix} 1&B_1+B_2\\0&1
\end{pmatrix}
$$
Also
$$
\begin{pmatrix}A&0\\0&A^{t-1}\end{pmatrix}
\begin{pmatrix} 1&B\\0&1
\end{pmatrix}
\begin{pmatrix}A&0\\0&A^{t-1}\end{pmatrix}^{-1}=
\begin{pmatrix} 1&ABA^t\\0&1
\end{pmatrix}
$$

Thus, we must prove, that
the there is no characters
$$\chi(B_1+B_2)=\chi(B_1)\chi(B_2)$$
on the additive group of
 symmetric integer matrices
$B$ such that
$$
\chi (ABA^t)=\chi(B)\qquad \text{for $A\in\SL(n,\Z)$}
$$
Any character of $N(\Z)$ has the form
$$
\chi(B)=\prod_{i\ge j} y_{ij}^{b_{ij}}, \qquad y_{ij}\in\C^*
$$

Since the group $\SL_n(\Z)$
contains all the even permutations of coordinates,
our
character has the form
$$\chi(B)= u^{\tr B}\cdot v^{\sum_{i>j}b_{ij}}$$
Next,
\begin{multline*}
\begin{pmatrix}
1&1&0\\0&1&0\\0&0&1
\end{pmatrix}
\begin{pmatrix}
a_{11}&a_{12}&a_{13}\\
a_{12}&a_{22}&a_{23}\\
a_{13}&a_{23}&a_{33}
\end{pmatrix}
\begin{pmatrix}
1&0&0\\1&1&0\\0&0&1
\end{pmatrix}
=\\=
\begin{pmatrix}
a_{11}+2 a_{12}+a_{22}&a_{12}+a_{22}&a_{13}+a_{23}\\
a_{12}+a_{22}&a_{22}&a_{23}\\
a_{13}+a_{23}&a_{23}&a_{33}
\end{pmatrix}
\end{multline*}
Hence, for all $a_{ij}\in\Z$ we have
$$
u^{2a_{12}+a_{22}} v^{a_{22}+a_{23}}=1
$$
Thus $u=v=1$. \hfill $\square$

\smallskip


{\bf \punct\label{ss4.5} Proof of Theorem \ref{th4.2}.}
Assume that we have some linearization
$S^\circ_\alpha$
of the representation
$S_\alpha$ on $\Sp(2n,\R)$.
By Lemma \ref{l4.4},
this linearization
 is rigidly defined on the matrices
$\begin{pmatrix}1&B\\0&1\end{pmatrix}$
by the formula
\begin{equation}
S_\alpha^\circ\begin{pmatrix}1&B\\0&1\end{pmatrix}
f(z)=f(z+B)
\label{eq4.6}
\end{equation}

Now let us consider the subgroup $\SL(2,\Z)=\Sp(2,\Z)\subset\Sp(2n,\Z)$
consisting of the
$[1+(n-1)+1+(n-1)]\times [1+(n-1)+1+(n-1)]$ matrices
$$
\begin{pmatrix}

a& 0&b&0\\
0&1&0&0\\
c&0&d&0\\
0&0&0&1
\end{pmatrix}
$$
To be short,
below we will write
$\begin{pmatrix}a&b\\c&d\end{pmatrix}$.
Let us show that
having  condition (\ref{eq4.6}),
 we can not trivialize the cocycle on
$\SL(2,\Z)$.

Denote
$$
I=\begin{pmatrix}0&-1\\1&0\end{pmatrix},\qquad
J=\begin{pmatrix}0&-1\\1&-1\end{pmatrix},\quad
K=\begin{pmatrix} 1&-1\\0&1\end{pmatrix}
$$
It can be easily be checked that
\begin{equation}
I^4=1,\quad J^3=1,\quad IK=J
\label{eq4.7}
\end{equation}

Our operators $S^\circ_\alpha(\cdot)$
have the form
\begin{align*}
S^\circ_\alpha(K)f(z)&=f(z-1),\\
S^\circ_\alpha(I)f(z)&=\theta\cdot f(-1/z)z^{-\alpha},\\
S^\circ_\alpha(J)f(z)&=\theta'\cdot f(-1-1/z) z^{-\alpha}
\end{align*}
where
$$
z^{-\alpha}=|z|^{-\alpha} \exp(-i\alpha \arg z);
\qquad 0<\arg z<\pi
$$
and $\theta$, $\theta'\in\C^*$ are some unknown constants.
The equation $IK=J$ implies $\theta'=\theta$.
Also
\begin{equation}
S^\circ_\alpha(I)^4=\theta^4 \exp(-2\alpha\pi i)
,\qquad
S^\circ_\alpha(J)^3=\theta^3\exp(-2\alpha\pi i)
\label{eq4.8}
\end{equation}
Since the both these operators equal 1,
$\theta=1$, $\exp(-2\alpha\pi i)=1$.
We obtaine a contradiction, since
 $\alpha\notin \Z$.

\smallskip

{\sc Remark.} An
 evaluation of powers in (\ref{eq4.8}) can be simplified in the following
 way. First, the point $i$ is a fix point of the transformation
 $z\mapsto -1/z$. Hence we can follow only
 the values of $S^\circ_\alpha(I)^k f(z)$ in this point.
 In the second case the fixed point is
 $\lambda=\exp(2\pi i/3)$.

{\sf Math.Phys. Group,
Institute of Theoretical and Experimental Physics,

B.Cheremushkinskaya, 25, Moscow 117259

\& University of Vienna, Math. Dept.,
Nordbergstrasse, 15, Vienna 1090, Austria}

neretin@mccme.ru

\end{document}